\documentclass[12 pt]{article}

\usepackage{amssymb} 
\usepackage{amsmath} 
\usepackage{latexsym}
\usepackage{graphicx}

\usepackage{enumitem}

\graphicspath{ {images/} }
\usepackage{mathtools}

\DeclarePairedDelimiter\floor{\lfloor}{\rfloor}

\providecommand{\keywords}[1]{\textbf{\textit{Keywords:}} #1}

\usepackage{tikz}
\usetikzlibrary{calc}
\usepackage{pgf,tikz}
\usetikzlibrary{arrows}
\usepackage{amscd}
\usepackage{relsize}

\usepackage{tikz-cd}

\makeatletter
\newcommand*{\rom}[1]{\expandafter\@slowromancap\romannumeral #1@}
\makeatother

\makeatletter
\def\blfootnote{\gdef\@thefnmark{}\@footnotetext}
\makeatother

\usepackage{ color, amsmath, amssymb, amsfonts, amscd, graphicx,latexsym, hyperref}
\usepackage{caption, subcaption}
\usepackage[pdf]{pstricks}
\usepackage{makeidx}

\usepackage[affil-it]{authblk}

\usepackage{amsthm}

\newtheorem{theorem}{Theorem}

\newtheorem{remark}{Remark}

\newtheorem{lemma}{Lemma}

\newtheorem{definition}{Definition}

\newtheorem{corollary}{Corollary}

\newcommand{\R}{\mathbb R}
\newcommand{\Z}{\mathbb Z}
\newcommand{\C}{\mathbb C}

\newcommand{\PP}{\frak P}

\title{Signature calculation of the area Hermitian form on some spaces of polygons}
\author{ \.{I}smail Sa\u{g}lam\thanks{Electronic address: \texttt{isaglamtrfr@gmail.com}  }}
\affil{Adana Alparslan Turkes Science and Technology University, \\ Institut de Recherche Mathematique Avanc\'ee, CNRS et Universit\'e de Strasbourg}

\date{}

\begin{document}

\maketitle

\begin{abstract}

This chapter is motivated by the paper by Thurston on triangulations of the sphere and singular flat metrics on the sphere. Thurston locally parametrized the moduli space of singular flat metrics on the sphere with prescribed positive curvature data by the complex hyperbolic space of appropriate dimension. This work can be considered as a generalization of signature calculation of the Hermitian form that he made in his paper.

The moduli space of singular flat metrics having unit area on the sphere with prescribed curvature data can be locally parametrized by certain spaces of polygons. This can be done by cutting singular flat spheres through length minimizing geodesics from a fixed singular point to the others. In that case the space of polygons is a complex vector space of dimension $n-1$ when there are $n+1$ singular points. Also there is natural area Hermitian form of signature $(1,n-2)$ on this vector space. In this chapter we calculate the signature of the area Hermitian form on some spaces of polygons which locally parametrize the moduli space of singular flat metrics having unit area on the sphere with one singular point of negative curvature. The formula we obtain depends only on the sum of the curvatures of the singular points having positive curvature.

This paper will appear in the book "In the tradition of Thurston, Vol. \rom{2}", Springer, 2022.
\end{abstract}
\keywords{singular flat metric, singular  flat surface, polygon, Hermitian form, Alexandrov Unfolding Process}
\blfootnote{\textup{} \textit{Mathematics Subject Classification}:
	51F99, 57M50, 15A69}

\section{Introduction}

Let $S$ be a compact singular flat surface\index{singular flat surface} perhaps with boundary. By this we mean that 

\begin{itemize}
	\item each interior point of $S$ has a neighborhood isometric to a neighborhood of the apex of a standard cone,
	\item
	near each boundary point $S$ possesses the geometry of a surface obtained by cutting a cone through a line passing from the apex. 
\end{itemize}

\noindent There is a well-defined notion of angle for each point $p\in S$. Let us denote the angle at $p$ by $\theta(p)$. If $p$ is an interior point we define the curvature at $p$ to be $\kappa_p=2\pi-\theta(p)$. If $p$ is a boundary point, then  $\kappa_p=\pi-\theta(p)$. Note that a point $p \in S$ is called singular if $\kappa_p\neq 0$. Otherwise $p$ is called non-singular. 

In \cite{WT} Thurston considered the moduli space of singular flat structures on the sphere 
with prescribed positive curvature data. He showed that if the number of singular points on the sphere is $n$, then the moduli space is a complex hyperbolic manifold of dimension $n-3$. To achieve this, he considered a singular flat metrics on the sphere and triangulated the sphere so that the vertices of the triangulation are exactly the singular points of the metric. Observing that nearby singular flat metrics admit the same combinatorial triangulation, he obtained local coordinates from the moduli space to the projectification 
of the positive part of a certain cocycle space 
equipped with the Hermitian form induced by the area of a singular flat structure. He showed that the signature of the Hermitian form is $(1,n-3)$, where $n$ is the number of singular points of the sphere. Note that this implies the local coordinates are from the moduli space to the complex hyperbolic space of dimension $n-3$. In this chapter we make a similar signature calculation for the case of singular flat spheres with one singular point of negative curvature. Now we return to the theory of triangulation of a singular flat surface $S$.

 It is well-known that a singular flat surface  $S$ can be triangulated by  Euclidean triangles. See \cite{MT1}. It follows that one can obtain any compact singular flat surface from a finite numbers of Euclidean triangles. However, the result in \cite{MT1} does not give us any constructive method to obtain a singular flat surface from triangles in the Euclidean plane.

There is a stronger result which says that a triangulation with a minimum number of triangles exists. More precisely, this theorem states that $S$ has a triangulation whose vertex set coincides with the set of singular points of the $S$. See \cite{ivan}, \cite{tahar} and \cite{ben-son} for proofs of this fact. However, even by this method, it is not clear  how one can construct a singular flat surface from a bunch of triangles in the Euclidean plane.

Another way to construct singular flat surfaces is to use flat disks\index{flat disk} instead of triangles. Note that a flat disk is a singular flat surface which is homeomorphic to a closed disk and has no singular interior points. It is not difficult to see 
that for any compact singular flat surface $S$, there exists a flat disk $D$ so that $S$ can be obtained by gluing some of the edges of $D$ appropriately. In \cite{ben-son}, it was shown that $D$ can be chosen so that it has a minimum number of edges. Equivalently, there exists a finite number of simple geodesic arcs on $S$ which can intersect only at their endpoints so that when we cut $S$ through these arcs we get a flat disk.

When $S$ has genus 0 and each point in $S$ has angle less than $2\pi$, there is a nice way to construct $S$ from a flat disk. Assume that $S$ has  $n+1$ singular points $P_0,P_1,\dots, P_{n}$ of positive curvature; that is, the angle at $P_i$, $\theta(P_i)$, is less than $2\pi$ for each $i$. Let $l_i$ be a length minimizing geodesic joining $P_0$ and $P_i$ for each $i>0$. Then it can be shown that each $l_i$ is simple and $l_i$ intersects $l_j$ only at $P_0$ when $i\neq j$. It follows that we get a flat disk if we cut $S$ through $l_1,l_2,\dots,l_n$. The Alexandrov Unfolding  Theorem\index{Alexandrov unfolding theorem} \cite{AD} states that this polygon can be embedded into the Euclidean plane. Therefore, any singular flat sphere with angle less than $2\pi$ at each  singular point can be obtained from a polygon in the Euclidean plane.

\begin{figure}
	
	\begin{center}
		\includegraphics[scale=0.7]{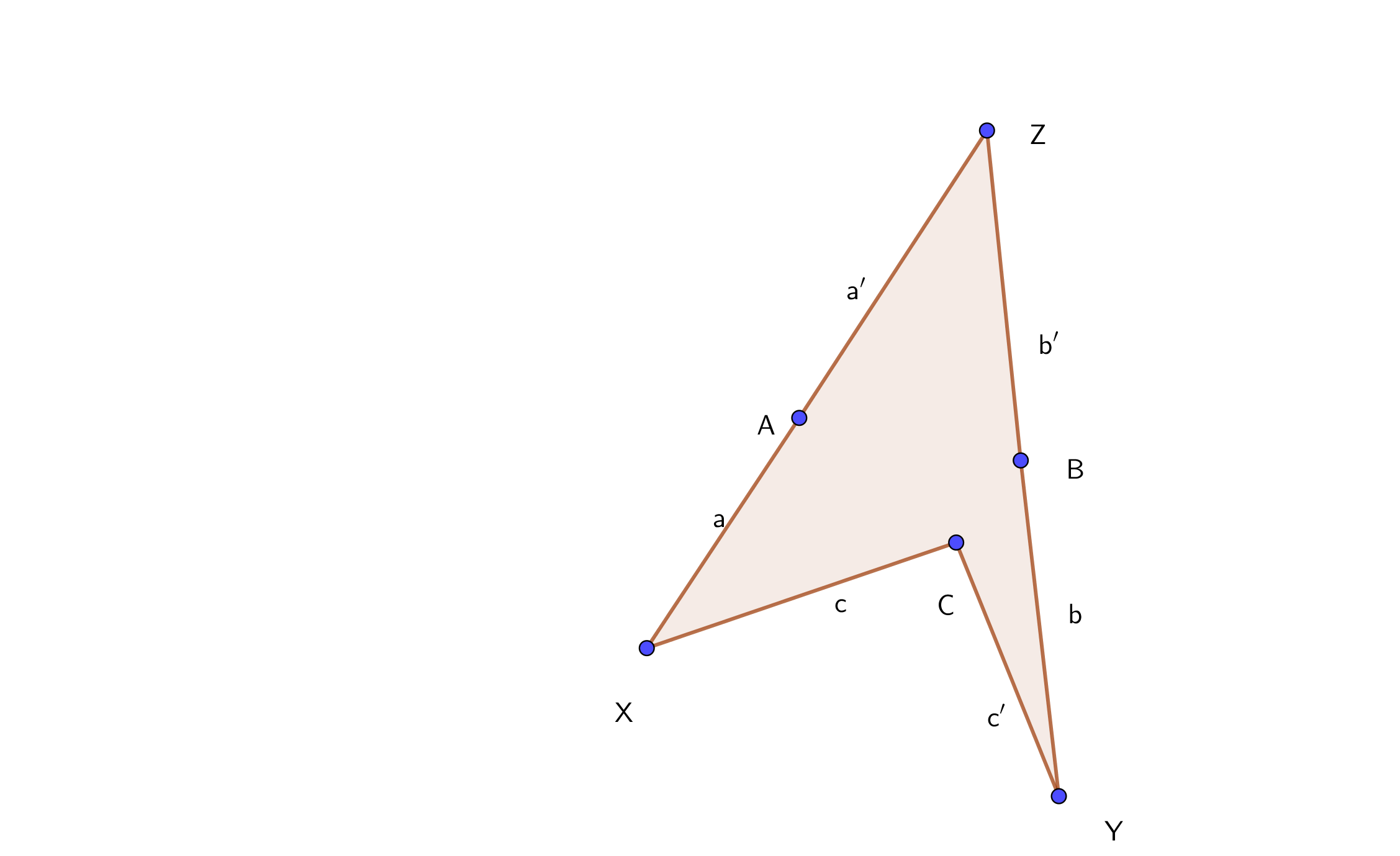}
	\end{center}
	
	\caption{A flat sphere with 4 singular points from a planer polygon. }
	\label{4tekil}
\end{figure}

For example, consider the polygon in Figure \ref{4tekil}. The edges that are denoted by the same letters have the same length. The polygon has two vertices with angle $\pi$: $A$ and $B$. Also, it has one vertex having angle $\frac{3\pi}{2}$. If we glue $a$ and $a'$, $b$ and $b'$, $c$ and $c'$, then we get a flat sphere. Note that this flat sphere has 4 singular points and 3 of them are obtained from $A$, $B$ and $C$. At these points the angles are $\pi$,$\pi$ and $\frac{3\pi}{2}$, respectively. Also, after this gluing operation the vertices $X$,$Y$ and $Z$ of the polygon come together to form a singular point of the flat sphere. This singular point has angle $\frac{\pi}{2}$. Note that the Alexandrov Unfolding Theorem states that any singular  flat sphere with 4 singular points of angle $\frac{\pi}{2},\frac{3\pi}{2},\pi,\pi$ can be obtained from such a polygon in the Euclidean plane.

Let us consider the general case again. That is, $S$ has genus 0 and $n+1$ singular points $P_0,\dots,P_n$ so that each $\theta(
	P_i)$ is less than $2\pi$. If we cut $S$ through $l_1,\dots,l_n$, then we get a flat disk which is isometric to a polygon. Note that this polygon is not arbitrary since it has certain properties. First of all, it has $2n$ vertices and $2n$ edges. It has $n$ vertices coming from $P_1,\dots, P_n$ so that the interior angles of the polygon at these vertices are   $\theta(P_1),\dots \theta(P_n)$. Let 

$$\kappa_i= 2\pi -\theta(P_i) \ \text{for}\ i=0,1,\dots,n.$$

\noindent The Gauss-Bonnet Formula \cite{MT3} implies that 

$$\sum_{i=0}^n\kappa_i= 4\pi$$

\noindent and we have

$$0<\kappa_0=4\pi-\sum_{i=1}^n\kappa_i<2\pi.$$ 

\noindent This polygon has $n$ vertices which are induced from $P_0$ and the sum of the angles at these vertices is equal to $P_0$.
 Furthermore, a vertex which comes from the singular point $P_i$ ( $1\leq i \leq n$) is incident to 2 edges of the same length. Let us move the polygon in the plane so that one of the vertices induced from $P_0$ is at the origin in the complex plane. If we assume that the vertices of the polygon having angles $\theta(P_1),\dots, \theta(P_n) $ are in counter-clockwise orientation, then these vertices give us an element $\hat{z}=(z_1,\dots,z_{2n})\in \C^{2n}$ so that 

\begin{enumerate}
	\item 
	$z_1=0$,
	\item
	$e^{\imath\kappa_k}(z_{2k-1}-z_{2k})= z_{2k+1}-z_{2k}, \ 1\leq k \leq n$.
\end{enumerate} 

\noindent Here, $\imath$ denotes the complex number $\sqrt{-1}$. The elements in $\C^{2n}$ which satisfy the above conditions form an $n-1$ dimensional complex vector space. We denote this vector space by $\frak{P}(\kappa)=\frak{P}(\kappa_1,\dots,\kappa_n)$. Note that the elements in $\frak{P}(\kappa) $ can be considered as possibly self-intersecting polygons. If such a polygon is positively oriented and does not have any self-intersection, then it represents a flat sphere. In this case the flat sphere and the polygon have the same area. If this polygon has coordinates $\hat{z}=(z_1,\dots,z_{2n})$, then its area is given by the following formula:

\begin{equation*}
\frac{\sqrt{-1}}{4}\sum_{i=1}^{2n-1}(z_i\overline{z}_{i+1}-z_{i+1}\overline{z}_i).
\end{equation*}

\noindent Therefore it is natural to consider the following area Hermitian form $h_A$ on $\frak{P}(\kappa)$:

$$h_A(\hat{z},\hat{w})= \frac{\sqrt{-1}}{4}\sum_{i=1}^{2n-1}(z_i\overline{w}_{i+1}-z_{i+1}\overline{w}_i).$$

This Hermitian form is called {\it area Hermitian form}\index{area Hermitian form}. Fillastre \cite{fillastre} computed the signature of the area Hermitian form and showed that it is $(1,n-2)$. Note that this computation is consistent with the one that was done by Thurston \cite{WT}. 

In this chapter, we calculate the signature of the area Hermition form $h_A$ on $\frak{P}(\kappa)=\frak{P}(\kappa_1,\dots,\kappa_n)$ by dropping some conditions on the curvature data $\kappa=(\kappa_1,\dots,\kappa_n)$. As before, we assume that $0<\kappa_i<2\pi$ for each $1\leq i \leq n$. But we do not require that

$$0<\kappa_0= 4\pi-\sum_{i=1}^{n}\kappa_i<2\pi.$$ 

\noindent That is, $\kappa_0$ may be any number less than $4\pi$. This has the following geometric significance. Assume that $\kappa_1+\dots+\kappa_n>4\pi$ and $\hat{z}\in \frak{P}(\kappa)$ is a positively oriented polygon. Then we can obtain a flat sphere by identifying equal edges of this polygon appropriately. This flat sphere has $n+1$ singular points and one of them has angle 

$$2\pi -\kappa_0=2\pi -(4\pi-(\kappa_1+\dots+\kappa_n))>2\pi.$$

\noindent Therefore, at least some part of the moduli   space of flat spheres with exactly one singular point having angle greater than $2\pi$ can be parametrized by using $\frak{P}(\kappa)$, for some $\kappa$. This gives us the hope to endow the moduli space of flat spheres with prescribed curvature data with new geometric structures.

The paper \cite{oshika} is closely related to the present one. Indeed our work can be considered as a  generalization of this paper by
Nishi annd Ohshika. In that paper the authors calculated the signature of the area Hermitian form 
on $\frak{P}(\kappa)$, where $\kappa=(\pi,\dots,\pi)$. This calculation led them to put a pseudo-metric on the moduli space/Teichm\"{u}ller space of flat metrics on the sphere with $n+1$ singular points, where the cone angles are $(n-2)\pi,\pi,\dots,\pi$. This new metric structure enabled them to put a pseudo-metric on the moduli space of hyperelliptic curves, since a hyperelliptic curve is a degree 2 branched cover of the complex projective line. We also point out that the papers \cite{star} and \cite{parametrisation} are closely related to the present one. Furthermore, ,in \cite{bavard-ghys}, the authors calculated the signature of a symmetric bilinear form on a space of polygons. The formulae and the proofs given in this chapter are similar. Finally, we note that Nishi \cite{nishi} addressed the question on the signature of a Hermitian form given
by the area function on the space of singular flat metrics on the sphere with conical singularities of possibly negative curvatures.

For more information about the geometry of polygons, see \cite{FRS} and \cite{Gro}. For more information about the geometry of flat surfaces, see \cite{MT1}, \cite{MT2}, \cite{Saglam} and \cite{MT3}. 
\section{Basic facts on Hermitian forms}
In this section we introduce the main facts that we use from the theory of the Hermitian forms\index{Hermitian form}. Let $V$ be a complex vector space. A Hermitian form on $V$ is a function

$$h:V\times V \to \C$$

\noindent such that 

\begin{enumerate}
	\item 
	$$h(\alpha u+ \beta v, w)=\alpha(u,w)+\beta h(u,w)$$
	\noindent for all $u,v,w \in V$ and for all $\alpha,\beta \in \C$,

	\item
	
	$$h(w,\alpha u+\beta v)=\bar{\alpha}h(w,u)+\bar{\beta}h(w,v)$$
	\noindent for all $u,v,w \in V$ and for all $\alpha, \beta \in \C$.
	
	\item
	
	$$h(u,v)=\overline{h(v,u)}$$
	\noindent for all $u,v \in V$.
\end{enumerate} 

Note that if $u\in V$, then $h(u,u)\in \R$. $h(u,u)$ is called the square-norm\index{square-norm} of $u$.

Assume that $V$ has dimension $n$ and $\mathcal{U}=\{u_1,\dots u_n\}$ is an ordered basis for $V$. Then the matrix 

$$H=[h(u_i,u_j)]$$

\noindent is called the matrix of $h$ in the ordered basis $\mathcal{U}$. This matrix has the property that
$H_{ji}=\overline{H}_{i j}$ for all $1\leq i,j \leq n$. That is, the transpose of the conjugate matrix of  $H$, which is denoted by $H^*$, is equal to $H$. Note that such a matrix is called Hermitian matrix\index{Hermitian matrix}.

\begin{definition}
	The rank of a Hermitian form\index{rank of a Hermitian form} $h$ is the rank of the matrix  $H$. It is denoted by $Rank(h)$.
\end{definition}

\begin{definition}
	A Hermitian form $h$ on an $n$-dimensional complex vector space $V$ is called non-singular if $Rank(h)=n$.
\end{definition}

Note that $h$ is non-singular if and only if for each $u \neq 0 \in V$ there exist $v \in V$ such that $h(u,v)\neq 0$. 
\begin{definition}
	Let  $h$ be a Hermitian form on a complex vector space $V$. 
	\begin{itemize}
		\item
	$h$ is called positive definite\index{positive definite Hermitian form} if $h(u,u)>0$ for all $u\neq 0 \in V $.
	\item
		$h$ is called negative definite\index{negative definite Hermitian form} if 
		$h(u,u)<0$ for all $u \neq 0 \in V$.
	\end{itemize}
\end{definition}
 The following theorem asserts that each finite-dimensional complex vector space has a basis $\mathcal{U}$ such that the matrix of $h$ in $\mathcal{U}$ is diagonal.
 
 \begin{theorem}
 	Let $V$ be a finite-dimensional complex vector space and let $h$ be a Hermitian form on $V$. Then there is an ordered basis for $V$ in which $h$ is represented by a diagonal matrix.
 \end{theorem}

\noindent One can sharpen the previous theorem so that the entries of the diagonal matrix are $1,-1$ or $0$. Here is the precise statement.

\begin{theorem}
	\label{main-signature}
	Let $V$ be an $n$-dimensional complex vector space and $h$ be a Hermitian form on $V$ which has rank $r$. Then there is an ordered basis $\{u_1,\dots,u_n\}$ for $V$ in which the matrix of $h$ is diagonal and such that
	
	$$h(u_j,u_j)=\pm 1  \ \text{for all}\ j=1,\dots r,$$ 
	
		$$h(u_j,u_j)= 0  \ \text{for all}\ j=r+1,\dots ,n.$$ 
		
		\noindent Furthermore, the number of vectors $u_j$ such that $h(u_j,u_j)=1$, $h(u_j,u_j)=-1$ and $h(u_j,u_j)=0$ are independent of the choice of basis.
	
\end{theorem}

We say that two elements $u,v \in V$ are orthogonal if $h(u,v)=0$.  Let $W$ be a subspace of $V$. Let us define $W^{\perp}$ as  the subspace of $V$ which consists of elements of $V$ that are orthogonal to each element in $W$. Note that $V^{\perp}$ has dimension $n-r$ and it has a basis $\{u_{r+1}\dots,u_n\}$. 

Let $U$ and $W$ be subspaces of $V$ such that

\begin{enumerate}
	\item 
	$U \cap W= 0$;
	\item 
	any element in $v\in V$ can be written as $v=u+w$, where $u\in U$ and $w \in W$;
	\item
	$h(u,w)=0$ for all $u \in U$ and $w \in W$.
\end{enumerate}
In this case we write 

$$V=U\bigoplus W,$$
\noindent and call $W$  an orthogonal complement\index{orthogonal complement} of $U$ in $V$.

Let us denote the cardinality of a set $A$ by  $\lvert A \rvert$. 

\begin{definition}
	Let $V$ be an $n$-dimensional complex vector space and $h$  a Hermitian form on $V$. Let $\mathcal{U}=\{u_1, \dots u_n\}$ be a basis as in Theorem \ref{main-signature}. We introduce the following quantities:
	
	\begin{enumerate}
		\item 
		$P(h)=\lvert \{u_j: h(u_j,u_j)=1 \}\rvert$;
		\item
		$N(h)= \lvert \{ u_j: h(u_j,u_j)=-1\}\rvert$;
		\item
		$Z(h)=\lvert \{ u_j : h(u_j,u_j)=0\}\rvert$.
	\end{enumerate}
 $(P(h),N(h))$ is called  signature\index{signature of a Hermitian form} of $h$. 
\end{definition}

Clearly, $P(h)+N(h)=Rank(h)$ and the dimension of $V^{\perp}$ is equal to $Z(h)$. Also, $P(h)+N(h)+Z(h)$ is equal to $n$. We have $Z(h)=0$ if and only if $h$ is non-singular. Furthermore, $h$ is positive definite if and only if $P(h)=n$, and $h$ is negative definite if and only if $N(h)=n$.

Now we define isomorphisms of complex vector space equipped with Hermitian forms. Let $V$ and $V'$ be complex vector spaces together with Hermitian forms $h$ and $h'$. We say that $V$ and $V'$ are isomorphic as vector spaces equipped with Hermitian forms if there is a vector space isomorphism $f: V \to V'$ such that

$$h(v,w)=h'(f(v),f(w))$$

\noindent for all $v,w \in V$. Note that $f$ is an isomorphism vector spaces having Hermitian forms if it satisfies the following weaker condition:

$$h(u,u)=h'(f(u),f(u))$$

\noindent for all $u \in V$. Note that the rank and signature  of a Hermitian form is invariant under isomorphisms. Also, two vector spaces having Hermitian forms are isomorphic if and only if they have the same signature and dimension.

\section{Spaces of polygons and signature calculation}
In this section, we introduce the spaces of polygons\index{space of polygons} that we consider. Each of these spaces is a complex vector space and admits  a natural {\it area Hermitian form} on it. We  calculate the signature of the area Hermitian form for each of these spaces.  
Let 
$$\kappa=(\kappa_1, \dots, \kappa_n), \ n>1,  \ 0 < \kappa_i < 2\pi, \ 1 \leq i \leq n,$$

\noindent be an $n-$tuple of real numbers. We will sometimes call it as curvature data. Let

\begin{align*}
\frak{P}(\kappa)=&\{\hat{z}=(z_1,\dots,z_{2n})\in \C^{2n}:\\
 &z_1=0, e^{\imath\kappa_i}(z_{2i-1}-z_{2i})=z_{2i+1}-z_{2i}, 1\leq i \leq n \}.
\end{align*}

\begin{figure}

	\begin{center}
	\includegraphics[scale=0.4]{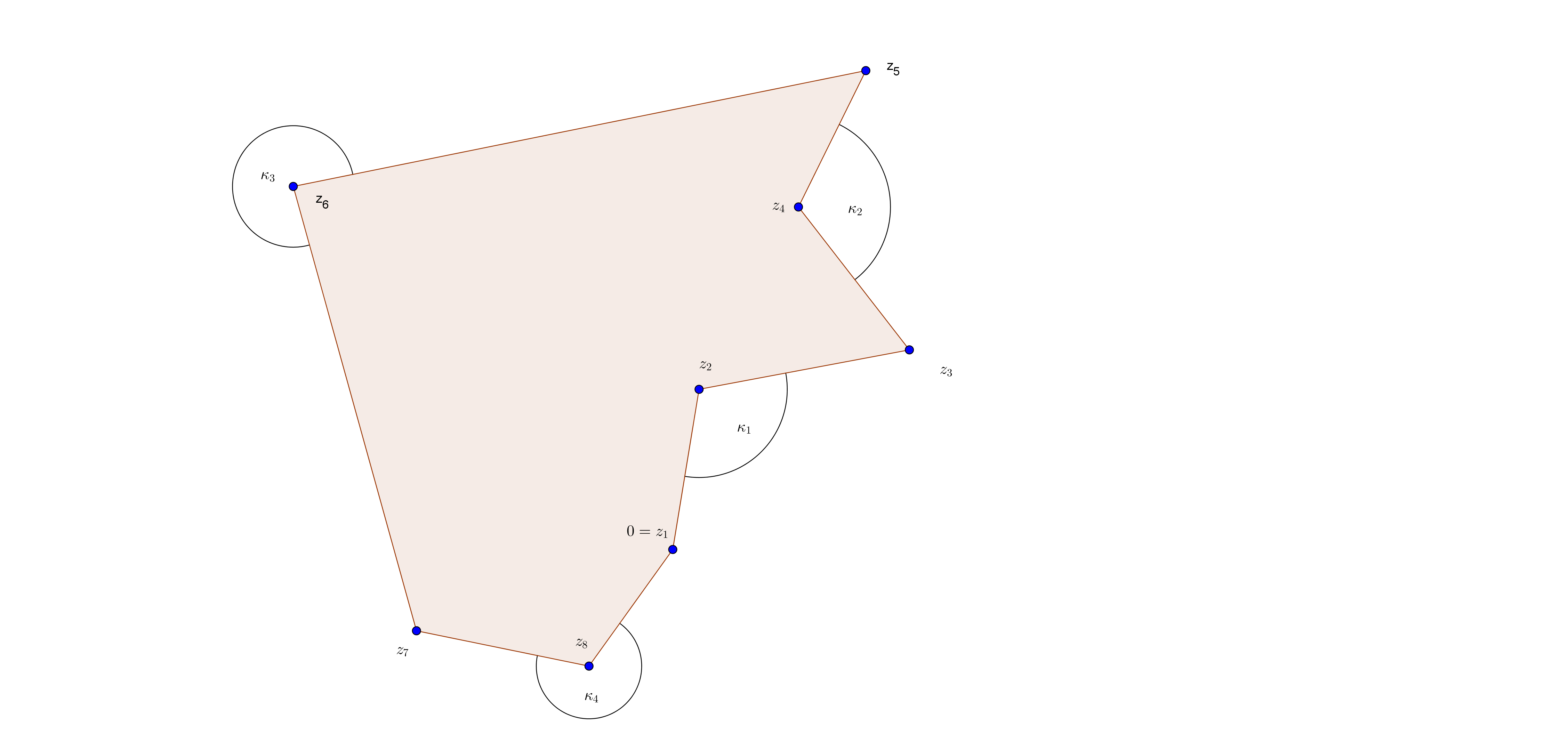}
	\end{center}

\caption{An element of  $\frak{P}(\kappa_1,\kappa_2,\kappa_3,\kappa_4)$ as a polygon in the complex plane. Note that  an edge of of  an element of $\frak{P}(\kappa_1,\kappa_2,\kappa_3,\kappa_4)$ considered as a polygon may intersect  another edge. }
	\label{basicpolygon}
\end{figure}

\noindent $\frak{P}(\kappa)$ can be thought as the set of oriented polygons 

$$z_1 \rightarrow z_2 \rightarrow z_3 \dots z_{2n}\rightarrow z_1.$$
Note that each element in $\PP(\kappa)$  has an outer angle $\kappa_i$ at the vertex $2i$, $z_{2i}$, where the outer angle is the angle between the vectors $z_{2i-1}-z_{2i}$ and $z_{2i+1}-z_{2i}$ measured  counter-clockwise. Also,  for all $1 \leq i \leq n$ and for all 
$z \in \frak{P}(\kappa)$, $\lvert z_{2i-1}-z_{2i}\rvert= \lvert z_{2i+1}-z_{2i}\lvert$. See Figure \ref{basicpolygon}.

\begin{remark}
	Dimension of $\frak{P}(\kappa)$ is $n-1$ since each element $\hat{z}=(0,z_2,\dots,z_{2n})$ is determined by its coordinates $z_3,\dots, z_{2n-1}$.
\end{remark}

\subsection{The area Hermitian form and the formula for its signature}
\label{signature}

In this section we introduce the Hermitian form that we are interested in. Also we give the formula for its signature.
 Consider the  Hermitian form 

\begin{align*}
h_A(\hat{z},\hat{w})=\frac{\sqrt{-1}}{4}\sum_{i=1}^{2n-1}(z_i\bar{w}_{i+1}-z_{i+1}\bar{w}_i)
\end{align*}
on $\frak{P}(\kappa)$, where $\hat{z}=(z_1,\dots,z_{2n}), \hat{w}=(w_1, \dots,w_{2n})$. We know that if $\hat{z}$ is  a simple polygon, then the area of $z$ is just the square-norm of $\hat{z}$, $h_A(\hat{z},\hat{z}) $. Therefore this form is called  area Hermitian form.

We now state the formula that we prove. First, we introduce some notation.

 Let  

\[
\epsilon(\kappa):=
\begin{cases}
1 &  \ \ \text{if}    \  \sum_{i=1}^n\kappa_i \in  2\pi \Z,\\
0 & \text{else.}
\end{cases}
\]

\noindent We denote the cardinality of a set $A$ by $\lvert A \rvert$. Let

\begin{align}
\label{defqkappa}
q(\kappa)=\big\lvert \{i: 1\leq i < n, \floor*{\sum_{k=1}^{i+1}\frac{\kappa_k}{2\pi}}=\floor*{\sum_{k=1}^i \frac{\kappa_k}{2\pi}}    \} \big\rvert,
\end{align}
and 
\begin{align}
\label{defpkappa}
p(\kappa)= n-1-q(\kappa)-\epsilon(\kappa),
\end{align}

\noindent where  $\floor*{}$ is the floor  function.

\begin{lemma}
	\label{floor-f}
	Let $f(i)=\floor*{\sum_{k=1}^i \frac{\kappa_k}{2\pi}} $. Then
	$$n-1-q(\kappa)=f(n).$$
	In particular, $q(\kappa)$ only depends on $\sum_{k=1}^n\kappa_k$.
	\begin{proof}
		
		Note that $f(i+1)-f(i)=0 \ \text{or} \ 1$. The definition of $q(\kappa)$ implies that 
		
$$n-1=q(\kappa)+\big\lvert \{i: 1\leq i < n, \floor*{\sum_{k=1}^{i+1}\frac{\kappa_k}{2\pi}}\neq\floor*{\sum_{k=1}^i \frac{\kappa_k}{2\pi}}    \} \big\rvert.$$

Therefore, it follows that 
		
$$n-1-q(\kappa)=\big\lvert \{i: 1\leq i < n, \floor*{\sum_{k=1}^{i+1}\frac{\kappa_k}{2\pi}}\neq\floor*{\sum_{k=1}^i \frac{\kappa_k}{2\pi}}    \} \big\rvert$$
		
	$$=\big\lvert \{i: 1\leq i <n, f(i+1)\neq f(i)                 \}\big\rvert$$	
		
			$$=\big\lvert \{i: 1\leq i <n, f(i+1)- f(i) =1                \}\big\rvert$$	
		$$=f(n)-f(n-1)+f(n-1)-f(n-2)+\dots f(2)-f(1)=f(n),$$
	\noindent since $f(1)$ is equal to 0. The particular case is obvious.
		
	\end{proof}
\end{lemma}

\begin{lemma}
	\label{permutationpq}
	If $\sigma$ is a permutation of $\{1,\dots,n\}$ and $\kappa(\sigma)=(\kappa_{\sigma(1)},\dots,\kappa_{\sigma(n)})$, then 
	
	$$q(\kappa)=q(\kappa({\sigma})), \ \textup{and}\ p(\kappa)=p(\kappa({\sigma}))\ \textup{and}\ \epsilon(\kappa)=\epsilon(\kappa(\sigma)).$$
	\begin{proof}
	It is clear from the definition of $\epsilon$ that $\epsilon(\kappa)=\epsilon(\kappa(\sigma))$. Also Lemma \ref{floor-f} implies that $q(\kappa)=q(\kappa(\sigma))$. Since 
	
	$$p(\kappa)=n-1-q(\kappa)-\epsilon(\kappa)$$	
	\noindent it follows that $p(\kappa)=p(\sigma(\kappa))$ for any permutation $\sigma$.
	\end{proof}
	
\end{lemma}

We will prove that the signature of the area Hermitian $h_A$ is 
          $$(p(\kappa),q(\kappa)),$$
\noindent that is, we will prove that

\begin{enumerate}
	\item 
	$P(h_A)=p(\kappa)=\floor*{\sum_{k=1}^n \frac{\kappa_k}{2\pi}}-\epsilon(\kappa)$;
	\item
	$N(h_A)=q(\kappa)= n-1-\floor*{\sum_{k=1}^n \frac{\kappa_k}{2\pi}} $;
	\item
	$Z(h_A)=\epsilon(\kappa)$.
\end{enumerate}

\noindent  Note that this will imply that 
$h_A$ is non-singular if and only if $\epsilon(\kappa)=0$.

\subsection{The case n=2}

\begin{figure} 
\centering
		\includegraphics[scale=0.4]{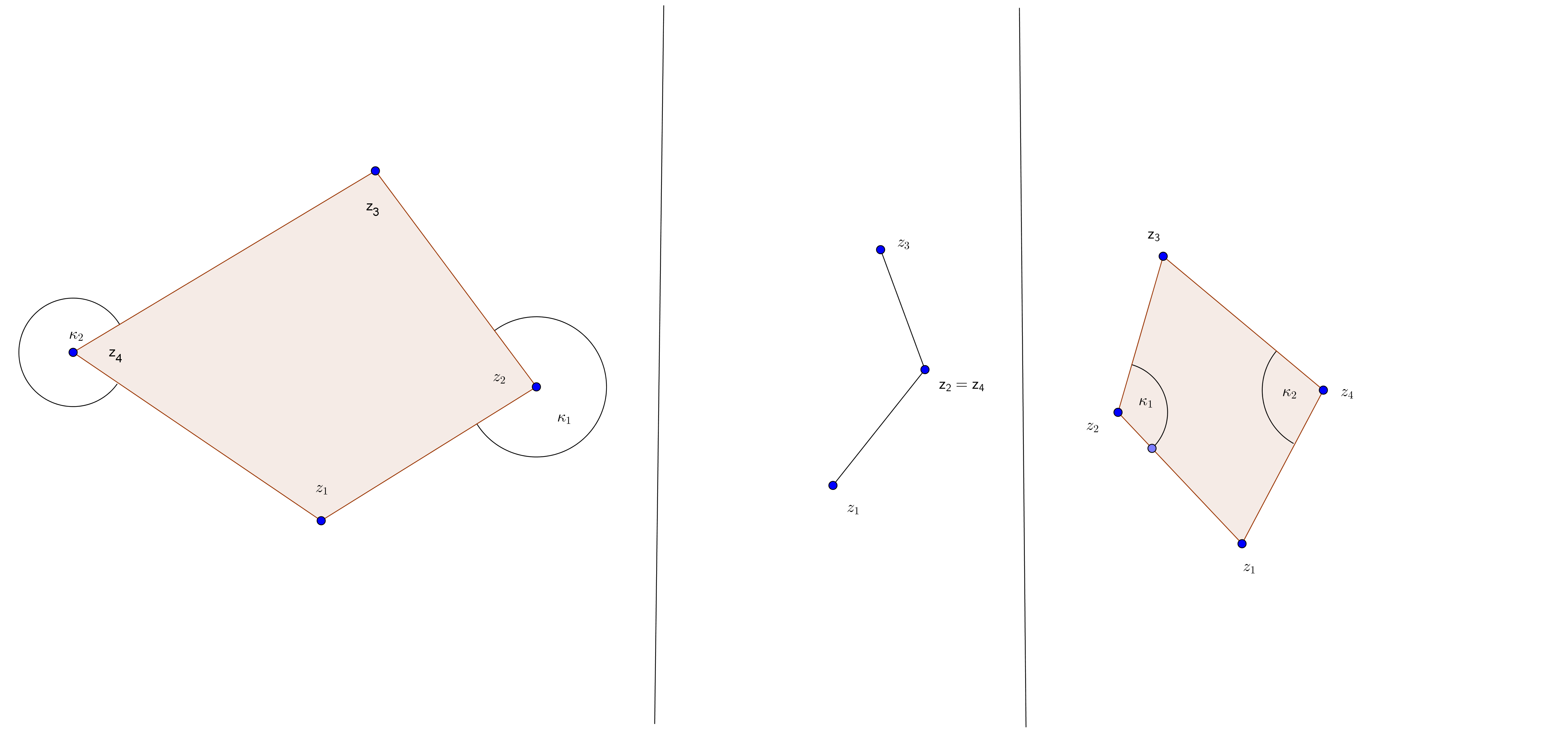}
		\caption{In each part of the figure, a generic element in $\PP(\kappa_1,\kappa_2)$ is given. In the leftmost picture, we have  $\kappa_1+\kappa_2>2\pi$. In the picture in the middle, we have $\kappa_1+\kappa_2=2\pi$. In the rightmost picture, we have $\kappa_1+\kappa_2<2\pi$. }
	\label{casesfor2}
\end{figure}

In this section we consider the case where $\kappa=(\kappa_1,\kappa_2)$. It follows that $\frak{P}(\kappa)$ has dimension 1. This means that the polygons that we consider have 4 edges and we can easily draw them. See Figure \ref{casesfor2}.
\begin{lemma}
	\label{signaturefor2}
If $n=2$, then the signature of the area Hermitian form  is 
$(p(\kappa), q(\kappa))$.
\begin{proof}

    There 
	are 3 cases to consider.
\begin{enumerate}\label{alph}
	\item 
	$\kappa_1+\kappa_2>2\pi$. From the left of  Figure \ref{casesfor2},  it is clear that the polygon corresponding to a non-zero element of $\frak{P}(\kappa)$ has  positive area. Since the area of the polygon corresponding to a non-zero element of $\frak{P}(\kappa)$ is the square-norm of that element, we see that each non-zero element of $\frak{P}(\kappa)$  has  positive square-norm. Since $\frak{P}(\kappa)$ is one-dimensional, it follows that the signature of the area  Hermitian form 
	is $(1,0)$. On the other hand, it follows directly from the definition of $p(\kappa)\  \text{and}\ q(\kappa)$ that  $(p(\kappa),q(\kappa))=(1,0)$. The result follows.
	\item 
	$\kappa_1+\kappa_2=2\pi$. In this case, every element in $\frak{P}(\kappa)$ has  area  0. See  the middle of  Figure  \ref{casesfor2} . Therefore, the signature is $(0,0)$. Also, it is clear that $p(\kappa)=q(\kappa)=0$.
	\item 
	$0<\kappa_1+\kappa_2<2\pi$.  In this case each polygon corresponding to a non-zero element is negatively oriented. This means that the orientation of these polygons is clockwise. Therefore, each non-zero element of this set has  negative area. See the right hand side of  Figure \ref{casesfor2}. 
	It follows that the signature of the form is $(0,1)$. Also, it is clear from the definition of $p(\kappa)$ and $q(\kappa)$  that $q(\kappa)=1$ and $p(\kappa)=0$.
\end{enumerate}
\end{proof}
\end{lemma}

\subsection{A special family of  polygons}
In this section, we assume that $\kappa=\kappa(n):=(\pi,\pi,\dots, \pi)$, where the curvature data $(\pi,\dots,\pi)$ has length $n$.
\begin{lemma}
	\label{allpi}
	Let $\kappa=\kappa(n)=(\pi,\dots,\pi)$.
	\begin{enumerate}
		\item
			If $n=2k+1$, then the signature of the area Hermitian form on $\PP(\kappa(n))$ is $(k,k)$.
			\item 
				If $n=2k$, then the signature of the area Hermitian form on $\PP(\kappa(n))$ is $(k-1,k-1)$.
					\end{enumerate}
	\begin{proof}
The case $n=2$ was proven in Lemma \ref{signaturefor2}. We will prove the lemma by induction on the length of curvature data $(\pi,\dots,\pi)$, $n$. 
Assume that $n\geq3$. We start with a useful observation. Consider the  map

\begin{align*}
\hat{z}=(0, z_2, z_3 \dots,z_{2n-1},z_{2n})\mapsto \hat{z}'=(0,z_{2n}, z_{2n-1},\dots, z_3,z_2)
\end{align*}
which sends $\PP(\kappa(n))$ to itself. It is clear that this map is a vector space isomorphism. It simply gives a polygon  the opposite orientation. Note that  the following formula holds:

$$h_A( \hat{z}',\hat{z}') = - h_A( \hat{z},\hat{z} ).$$ 

\noindent It follows that $P(h_A)=N(h_A)$. Consider the following vector subspace of $\PP(\kappa(n))$: 

$$\PP'=\{\hat{z} \in \PP(\kappa(n)): z_2=0 \}.$$

\noindent Note that $\PP'$ is an $(n-2)$-dimensional vector subspace. Also, consider the restriction of the area Hermition form on $\PP'$. It is easy to see that $\PP'$ and
$\PP(\kappa(n-1))$ are isomorphic as complex vector spaces with Hermitian forms. 

Consider a basis $\{u_1,\dots,u_r,u_{r+1},\dots,u_{r+s}, \dots,u_{n-2}\}$ for $\PP'$ so that

\begin{enumerate}
	\item 
	$h_A(u_i,u_j)=0$ if $i\neq j$;
	\item
	$h_A(u_i,u_i)=1$ if $1\leq i \leq r$;
	\item
	$h_A(u_i,u_i)=-1$ if $r+1\leq i \leq r+s$;
	
\item 

and $h_A(u_i,u_i)=0$ if $i>r+s$.
\end{enumerate}

Let 

\begin{enumerate}
	\item 
	$\PP'^+$ be the vector subspace spanned by $\{u_1,\dots,u_r\}$;
	\item
	$\PP'^-$ be the vector subspace spanned by $\{u_{r+1},\dots,u_{r+s}\}$; 
	\item
	$\PP'^{\perp}$ be the vector subspace spanned by $\{u_{r+s+1},\dots, u_{n-2}\}$.
\end{enumerate}

Assume that $n=2k$, $k\geq 2$. Then $\PP'$ is $(2k-2)$-dimensional and the induction hypothesis implies that the signature of the area Hermitian form on $\PP'$ is $(k-1,k-1)$.
Therefore 
$$\PP'=\PP'^+\bigoplus\PP'^-,$$
where $\PP'^+$ and $\PP'^-$ have dimension $k-1$. Consider an element $u\in \PP(\kappa(n))\setminus \PP'$. Applying the Gram\textendash Schmidt orthogonalization process if necessary, we can choose $u$ so that it is orthogonal to $\PP'$. Therefore the signature of the area Hermitian form on $\PP(\kappa(n))$ is 
\begin{enumerate}
	\item 
	$(k-1,k-1)$,
	\item
	$(k,k-1)$ or
	\item
	$(k-1,k)$.
\end{enumerate}

\noindent Since $P(h_A)=N(h_A)$, it follows that this signature is $(k-1,k-1)$.

Assume that $n=2k+1$. Then $\PP'$ is a $2k-1$ dimensional subspace of $\PP(\kappa(n))$. The induction hypothesis implies that the signature of the area Hermitian form on $\PP'$ is $(k-1,k-1)$. Therefore the dimension of $\PP'^+$ is $k-1$, the dimension of $\PP'^-$ is $k-1$ and the dimension of $\PP'^{\perp}$ is $1$. 

Take an element $v\in \PP(\kappa(n))\setminus \PP'$ so that $v$ is orthogonal to $\PP'^+$ and $\PP'^-$. We can assure this by the Gram\textendash Schmidt orthogonalization process. Let $W$ be the vector space spanned by $v$ and $\PP'^{\perp}$. It is clear that 
\begin{equation}
	\label{decomp}
\PP(\kappa(n))= \PP'^+\bigoplus\PP'^-\bigoplus W. 
\end{equation}

\noindent Assume that $W$ is orthogonal to $\PP(\kappa(n))$, or equivalently, that the area Hermitian form restricted to $W$ is zero. Let

$$w=(0,0,0,w_4,\dots,w_{2n})$$
\noindent be a generator of $\PP'^{\perp}$. There is an integer $l$ such that $w_k=0$ for $k<2l$ and $w_{2l}\neq 0$. Note that $2 w_{2l}=w_{2l+1}$. Consider the following element of $\PP(\kappa(n))$:

$$a=(0,\dots,0,1,2,1,0,\dots,0),$$

\noindent where the first 1 is in the coordinate $2l-2$. Then 

$$h_A(a,w)=\sqrt{-1}\bar{w}_{2l}\neq 0.$$

\noindent This is a contradiction. It follows that the area Hermitian form restricted to $W$ is not trivial. Therefore the signature of the area Hermitian form restricted to $W$ is

\begin{enumerate}
	\item 
	$(0,1)$,
	\item
	$(1,0)$ or
	\item
	$(1,1)$.
\end{enumerate}

Regarding the decomposition \ref{decomp}, it follows that the signature of $h_A$ on $\PP(\kappa(n))$ is

\begin{enumerate}
	\item 
	$(k-1,k)$,
	\item
	$(k,k-1)$ or
	\item
	$(k,k)$. 
\end{enumerate}
Since $P(h_A)=N(h_A)$, the signature is $(k,k)$.
	\end{proof}

\end{lemma}
The following corollary is an immediate application of Lemma \ref{allpi}.

\begin{corollary}
	\label{signatureallpi}
	If $\kappa=(\pi,\dots,\pi)$, then the signature of the area Hermitian form is $(p(\kappa),q(\kappa))$.
	\begin{proof}
		It is not difficult to see that
	\begin{itemize}
		\item 
		$(p(\kappa),q(\kappa))=(k,k)$ if $n=2k$, and
		\item
		$(p(\kappa),q(\kappa))=(k-1,k-1)$ if $n=2k$.
		
	\end{itemize}
Therefore the statement follows from Lemma \ref{allpi}.
	\end{proof}
\end{corollary}

\subsection{Cutting-Gluing Operations}
In this section, we explain why for any permutation $\sigma\in S_n$,
$\PP(\kappa)$ and $\PP(\kappa(\sigma))$ are isomorphic as complex vector spaces equipped with Hermitian forms. 

We prove the claim by using some cutting and gluing operations\index{cutting and gluing operation}. Before proceeding to the general case, we first introduce cutting and gluing operations on $\PP(\kappa)$, where $\kappa=(\pi,\pi,\pi)$. We know that $\PP(\kappa)$ is 2-dimensional and that the signature of $h_A$ on it  is $(1,1)$. See Lemma \ref{allpi}. Take an element $\hat{z}=(z_1,z_2,z_3,z_4,z_5,z_6) \in \PP(\kappa)$. Then 

$$\hat{z}=(z_1,z_2,z_3,z_4,z_5,z_6)= (0,\frac{z_3}{2},z_3,\frac{z_3+z_5}{2},z_5,
\frac{z_5}{2})$$

\noindent and 

$$h_A(\hat{z},\hat{w})=\frac{\sqrt{-1}}{4}(z_3\bar{w}_5-z_5\bar{w}_3).$$

\noindent By abusing notation, we denote an element $\hat{z} \in \PP(\kappa)$ as $\hat{z}=[[z_3,z_5]]$. Now consider a positively oriented element $\hat{z} \in \PP(\kappa)$. This element has positive square-norm, that is,

$$h_A(\hat{z},\hat{z})=\frac{\sqrt{-1}}{4}(z_3\bar{z}_5-z_5\bar{z}_3)>0.$$

\noindent Recall that $h_A(\hat{z},\hat{z})$ is the area of the corresponding polygon which actually is a triangle. 

\begin{figure}
	
	\begin{center}
		\includegraphics[scale=0.65]{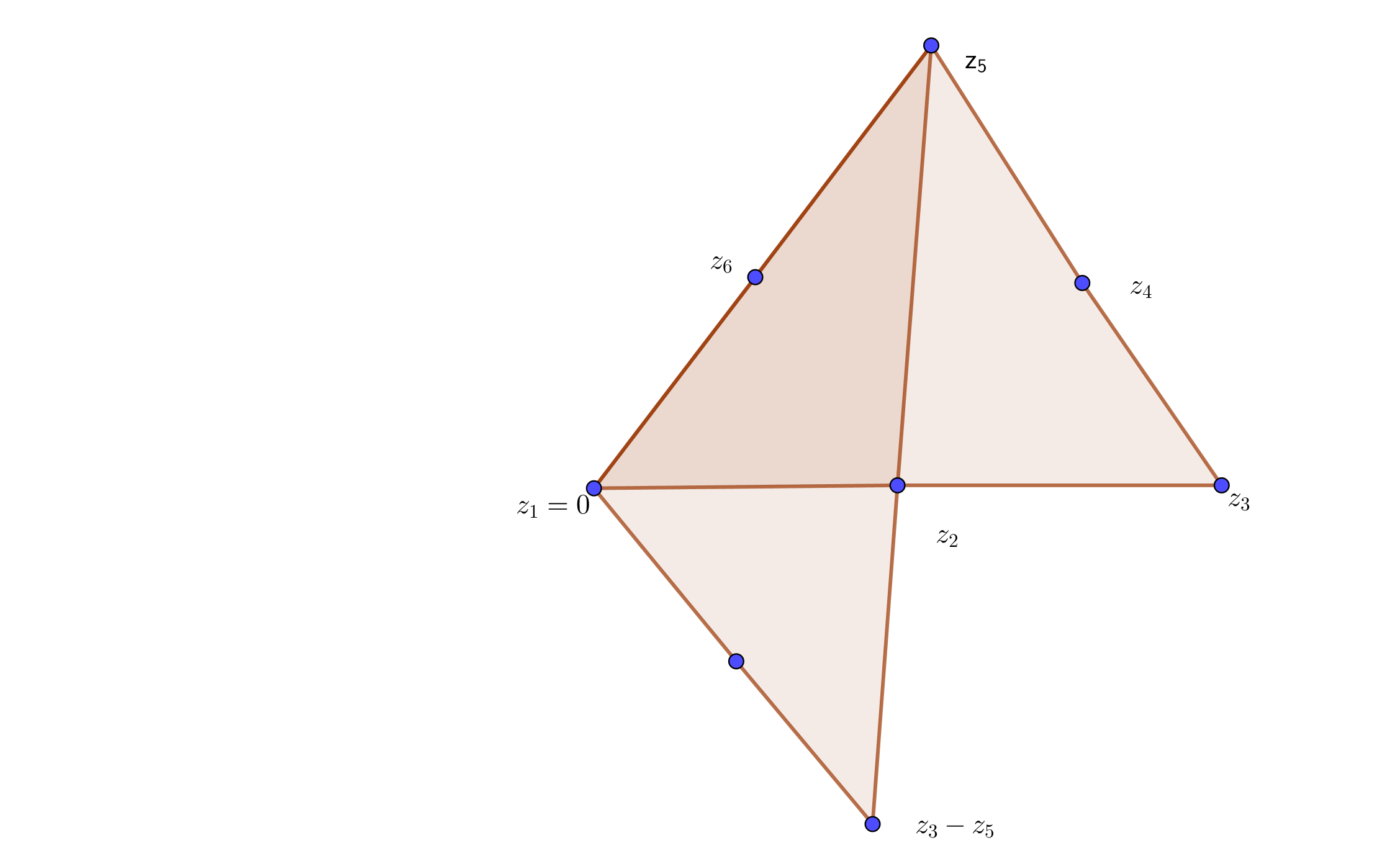}
	\end{center}
	
	\caption{Cutiing-gluing operation on $\PP(\pi,\pi,\pi)$. }
	\label{klasik}
\end{figure}

In Figure \ref{klasik}, consider the line segment $[z_2,z_5]$ and cut the triangle $[0,z_3,z_5]$ through this line segment to get two triangles. Glue the edges $[z_2,z_3]$ and $[0,z_2]$ by a rotation of angle $\pi$ around $z_2$. In this way, we get another element in $\PP(\kappa)$ having the same area with coordinates $[[z_3-z_5,z_5]]$. Therefore we have a map 

$$\PP(\kappa)\to \PP(\kappa)$$

\noindent sending $[[z_3,z_5]] \mapsto [[z_3-z_5,z_5]]$. Clearly this map is a vector space isomorphism and we realized that it respects the area Hermitian form.

\begin{figure}
	
	\begin{center}
		\includegraphics[scale=0.55]{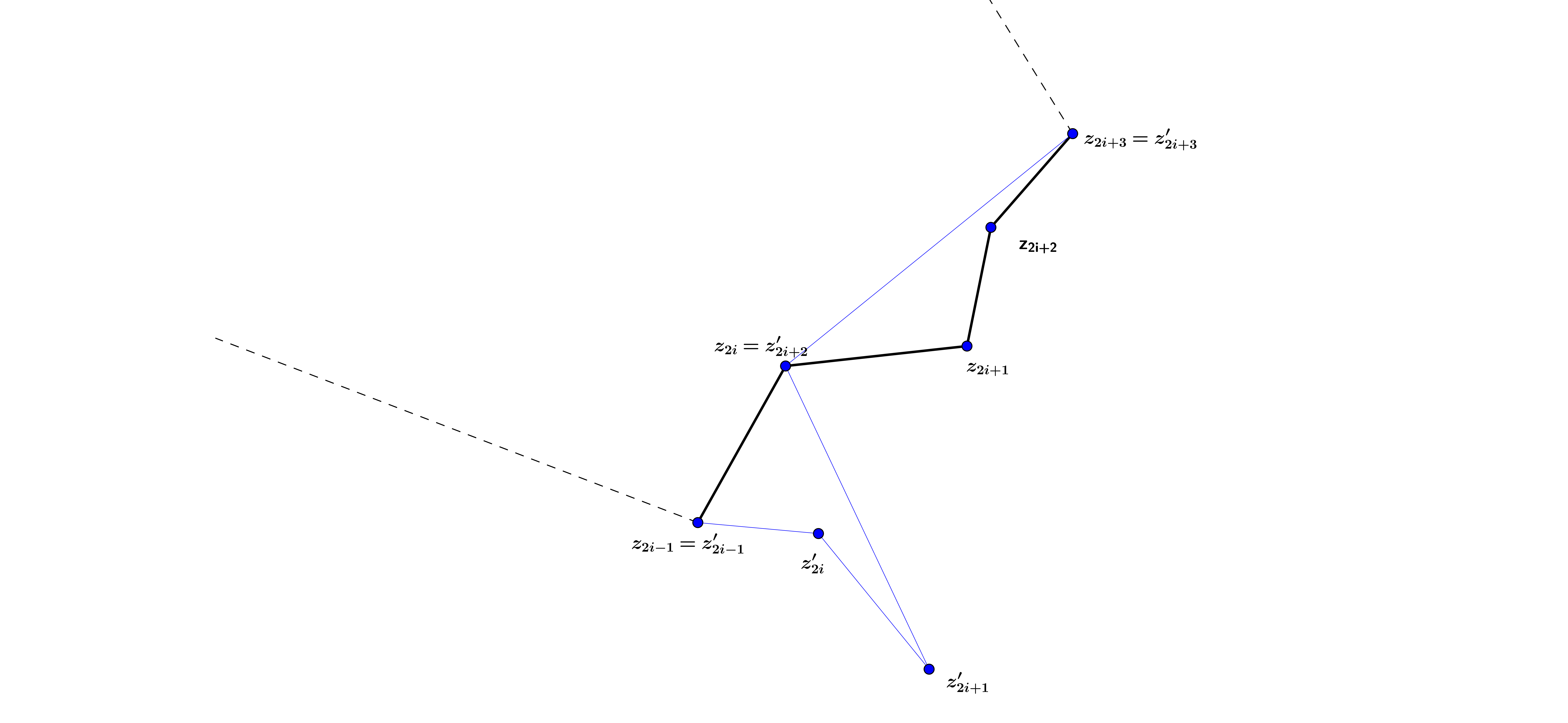}
	\end{center}
	
	\caption{We obtain an element  $\PP(\kappa(\sigma))$ from an element of  $\PP(\kappa).$ }
	\label{cutandglue}
\end{figure}

In general, even if the entries of the curvature data $\kappa$ are not equal, we can use these cutting-gluing operations to find isomorphisms between $\PP(\kappa)$ and $\PP(\kappa(\sigma))$.
 Note that it is enough to consider the cases for which $\sigma=(i,i+1)$ is a transposition to prove that $\PP(\kappa)$ and $\PP(\kappa(\sigma))$ are isomorphic. Take an element  $\hat{z}\in \PP(\kappa)$, and consider it as a polygon in the complex plane. Assume that the line segment joining $z_{2i+3}$ and $z_{2i}$ does not intersect the polygon except at its endpoints. Cut the polygon through the line segment joining $z_{2i+3}$ and $z_{2i}$. Glue the edge  $[z_{2i},z_{2i+1}]$ in the resulting quadrangle  with the edge  $[z_{2i},z_{2i-1}]$ of the polygon  to get the element $\hat{z}'$. In this way, we get an area-preserving map from a subset of $\PP(\kappa)$ to $\PP(\kappa(\sigma))$, where $\sigma=(i,i+1)$. Note that this map extends to an area-preserving linear map between 
 $\PP(\kappa)$ and $\PP(\kappa(\sigma))$. Indeed, this linear map is an isomorphism; one can reverse the cutting and gluing operations to get an inverse for the map.
See Figure \ref{cutandglue}.

Therefore, we have proved the following lemma.

\begin{lemma}
	\label{permutation-signature}
	\begin{enumerate}
		\item
		For any $\sigma \in S_n$, $\PP(\kappa)$ and $\PP(\kappa(\sigma))$ are isomorphic as complex  vector spaces equipped with Hermitian forms.
		\item
		The signature of the area Hermitian form on $\PP(\kappa)$ is equal to	the signature of the area Hermitian form on $\PP(\kappa(\sigma))$.
		\end{enumerate}
\end{lemma}
Note that these cutting-gluing operations were introduced in \cite{iso}.
\subsection{Signature calculation}

In this section we prove the signature formula. Let $\kappa=(\kappa_1,\dots,\kappa_n)$. Assume that $n>2$ and $\kappa_1+\kappa_2 < 2\pi$. Let

\[
\kappa_{12} =
\begin{cases}
	\kappa_1+\kappa_2 & \text{if $\kappa_1+\kappa_2<2\pi$} \\
	\kappa_1+\kappa_2-2\pi & \text{if $\kappa_1+\kappa_2>2\pi$} 
\end{cases}
\]

Also let $\kappa'=(\kappa_{12},\kappa_3, \dots, \kappa_n)$.
Consider the following $n-2$ dimensional subspace of $\PP(\kappa)$:

\begin{align*}
\overline{\PP}(\kappa)=\{\hat{z} \in \PP(\kappa): z_2=z_4 \}.
\end{align*}

\begin{figure}
	
	\begin{center}
		\includegraphics[scale=0.60]{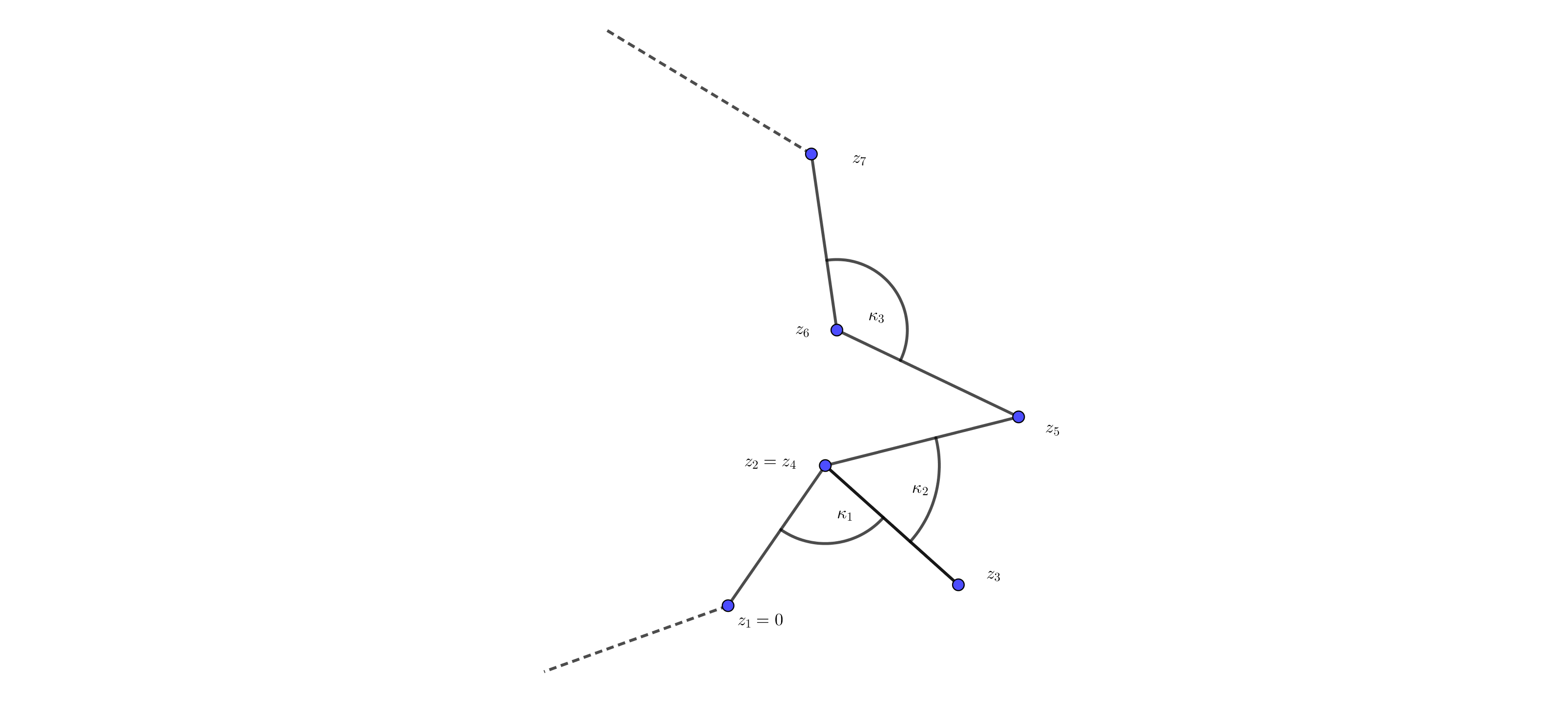}
	\end{center}
	
	\caption{A generic element of $\overline{\PP}(\kappa)$ when $\kappa_1+\kappa_2<2\pi$.}
		
	\label{generic}
\end{figure}

A generic element of $\overline{\PP}(\kappa)$ is shown in Figure \ref{generic}. It is not difficult to see that $\overline{\PP}(\kappa)$ together with the induced Hermitian form and $\PP(\kappa')$ are isomorphic. We want to find an orthogonal complement of $\overline{\PP}(\kappa)$ in $\PP(\kappa)$. In $\PP(\kappa)$, there is a unique element of the form

\begin{align*}
X=(0, -1, -1+e^{i\kappa_1},x,0,\dots,0).
\end{align*}
Since the angle at the fourth vertex, $x$, is $\kappa_2\neq 2\pi-\kappa_1$, it follows that $X \notin \bar{\PP}(\kappa)$ and

\begin{align*}
(-1+e^{\imath \kappa_1}-x)e^{\imath \kappa_2}=-x,\\
e^{\imath(\kappa_1+\kappa_2)}-e^{\imath\kappa_2}=x(e^{\imath\kappa_2}-1),\\
1-e^{-\imath\kappa_1}+x(-e^{-\imath\kappa_1}+e^{-\imath(\kappa_1+\kappa_2)})=0.
\end{align*}

\noindent On the other hand, any element in $\overline{\PP}(\kappa)$ is a constant multiple of an element of the form 

\begin{align*}
Y=(0,-1,-1+e^{\imath\kappa_1},-1,-1+e^{\imath(\kappa_1+\kappa_2)},\dots).
\end{align*}
It follows that 
\begin{align*}
 h_A(X,Y)=1-e^{-\imath\kappa_1}+x(-e^{-\imath\kappa_1}+e^{-\imath(\kappa_1+\kappa_2)})=0.
\end{align*}
Let $\C X$ denote the vector space generated by $X$. Clearly $\C X\equiv\PP(\kappa_1, \kappa_2) $. Note that  we have proved the following lemma.

\begin{lemma} If $n>2$ and $\kappa_1+\kappa_2\neq2\pi$, then
	\label{decomposition}
	$\C X \oplus \overline{\PP}(\kappa)=\PP(\kappa)$.
\end{lemma}

\begin{theorem}
The signature of the area Hermitian form for $\PP(\kappa)$ is
$(p(\kappa),q(\kappa))$. Also we have the following formulas for $p(\kappa)$ and $q(\kappa)$:

$$q(\kappa)=n-1-\floor*{\sum_{k=1}^n \frac{\kappa_k}{2\pi}}$$
\noindent and

$$p(\kappa)=\floor*{\sum_{k=1}^n \frac{\kappa_k}{2\pi}}-\epsilon(\kappa),$$

\noindent where 

\[
\epsilon(\kappa):=
\begin{cases}
	1 &  \ \ \text{if}    \  \sum_{i=1}^n\kappa_i \in  2\pi \Z,\\
	0 & \text{else.}
\end{cases}
\]
\begin{proof}
We prove the first part of the theorem by induction on $n$. If $n=2$ we know that the statement is true. See Lemma \ref{casesfor2}. Assume that $n>2$.
If all of the $\kappa_i$ are $\pi$, then we know the theorem is true. See 
Corollary \ref{signatureallpi}. Assume that not all of the $\kappa_i$ are $\pi$. It follows that there are $i,j \in \{1,2, \dots n\}$, $i\neq j$, such that $\kappa_i+\kappa_j\neq2\pi$. By Lemma
\ref{permutationpq} and Lemma \ref{permutation-signature}, we may assume that $\kappa_1+\kappa_2 \neq 2\pi$.

Assume that $\kappa_1+\kappa_2<2\pi$.
 Note that Lemma \ref{decomposition}  and the induction hypothesis on
the  vector spaces $\C X\equiv \PP(\kappa_1, \kappa_2)$,  $\overline{\PP}(\kappa)\equiv\PP(\kappa_1+\kappa_2,\kappa_3,\dots \kappa_n)$  imply that

\begin{align*}
N(h_A)&=q(\kappa_1,\kappa_2)+q(\kappa_1+\kappa_2,\kappa_3,\dots,\kappa_n)\\
&=\big\lvert \{i: 1\leq i < 2, \floor*{\sum_{k=1}^{i+1}\frac{\kappa_k}{2\pi}}=\floor*{\sum_{k=1}^i \frac{\kappa_k}{2\pi}} \} \big\lvert\\
&+\big\lvert\{i: 2\leq i < n, \floor*{\sum_{k=1}^{i+1}\frac{\kappa_k}{2\pi}}=\floor*{\sum_{k=1}^i \frac{\kappa_k}{2\pi}}      \}\big\rvert\\
&=\big\lvert \{i: 1\leq i < n, \floor*{\sum_{k=1}^{i+1}\frac{\kappa_k}{2\pi}}=\floor*{\sum_{k=1}^i \frac{\kappa_k}{2\pi}}    \} \big\rvert\\
&=q(\kappa).
\end{align*}
And we also have
\begin{align*}
P(\kappa)&=p(\kappa_1,\kappa_2)+p(\kappa_1+\kappa_2,\kappa_3,\dots\kappa_n)\\
&=1-q(\kappa_1,\kappa_2)-\epsilon(\kappa_1,\kappa_2)+n-2 \\
&-q(\kappa_1+\kappa_2,\kappa_3,\dots,\kappa_n)-\epsilon(\kappa_1+\kappa_2,\kappa_3,\dots,\kappa_n)\\
&=n-1-q(\kappa)-\epsilon(\kappa_1+\kappa_2,\kappa_3,\dots,\kappa_n)\\
&=n-1-q(\kappa)-\epsilon(\kappa)\\
&=p(\kappa).
\end{align*}

\noindent Therefore $(P(h_A),N(h_A))=(p(\kappa),q(\kappa))$. Note that a similar calculation holds for the  case $\kappa_1+\kappa_2>2\pi$.

The formulae for $p(\kappa)$ and $q(\kappa)$  follow easily from Lemma \ref{floor-f}.
\end{proof}
	\end{theorem}

\begin{corollary}Let $\kappa=(\kappa_1,\dots,\kappa_n)$.
	\begin{enumerate}
		\item 
		If  $2\pi<\sum_{i=1}^n\kappa_i<4\pi$, then the signature is $(1,n-2)$.
		\item
			If  $\sum_{i=1}^n\kappa_i=2\pi$, then the signature is $(0,n-2)$.
		\item
			If  $\sum_{i=1}^n\kappa_i=2\pi(n-1)$, then the signature is $(n-2,0)$.
		\item
			If $\sum_{i=1}^n\kappa_i<2\pi$, then the signature is $(0,n-1)$.
		\item
		If  $2\pi(n-1)<\sum_{i=1}^n\kappa_i<2\pi n$, then the signature is $(n-1,0)$.
			\item
		If  $2\pi(n-2)<\sum_{i=1}^n\kappa_i<2\pi(n-1)$, then the signature is $(n-2,1)$.
	\end{enumerate}
\end{corollary}

\section*{Acknowledgments}
The author is grateful to Susumu Tanabe, Muhammed Uluda\u{g}, Athanase Papadopoulos, Ke'nishi Ohshika and Haruko Nishi for their suggestions and help. 
He is financially supproted by T\"{U}B\.{I}TAK.


\end{document}